\documentclass[12pt,a4paper]{article}

\usepackage{latexsym,amssymb,amsthm}
\usepackage{url,graphics}
\usepackage[dvips]{graphicx,epsfig}

\title{\bf Bounds of a  number of leaves of spanning trees}
\author{ A.\,V.\,Bankevich  \and  D.\,V.\,Karpov }
\date{}

\begin{document}
\maketitle
\righthyphenmin=2
\renewcommand*{\proofname}{\bf Proof}
\newtheorem{thm}{Theorem}
\newtheorem{lem}{Lemma}
\newtheorem{cor}{Corollary}
\theoremstyle{definition}
\newtheorem{defin}{Definition}
\theoremstyle{remark}
\newtheorem{rem}{\bf Remark}

\def\N{{\rm N}}
\def\q#1.{{\bf #1.}}
\def\I{{\rm Int}}
\def\R{{\rm Bound}}
\def\mmin{\mathop{\rm min}}

\section{\bf Introduction. Basic notations}

We use standart notations. For a graph $G$ we denote the set of its vertices by  $V(G)$ and the set of its edges by $E(G)$. We use notations  $v(G)$ and $e(G)$ for the number of vertices and edges of $G$, respectively. 
All graphs in our paper do not contain loops  and multiple edges.

We denote the {\it degree} of a vertex $x$ in the graph $G$ by $d_G(x)$, and, as usual, denote the minimal vertex degree of the graph  $G$ by $\delta(G)$. 

Let $\N_G(w)$ denote the {\it neighborhood}  of a vertex  $w\in V(G)$ (i.e. the set of all vertices, adjacent to~$w$). We denote the girth of a graph $G$ (i.e. the length of its minimal cycle) by $g(G)$.

\goodbreak
\begin{defin}
For any connected graph~$G$ we denote by~$u(G)$  the maximal number of leaves in a spanning tree of the graph~$G$. 
\end{defin}

\begin{rem}
Obviously, if  $F$ is a tree, then   $u(F)$ is the number of its leaves.
\end{rem}

Several papers about lower bounds of $u(G)$ are published.  
In 1981 Storer~\cite{Storer} supposed, that  $u(G)>{1\over 4}v(G)$ as  $\delta(G)\ge 3$. 
Linial formulated a more strong conjecture:  $u(G)\ge {\delta(G)-2\over \delta(G)+1}v(G) + c$ as $\delta(G)\ge 3$, where a constant $c>0$ depends only on $\delta(G)$. This conjecture is reasoned by the fact that for every  $d\ge 3$ one can easily construct an infinite series of graphs with   minimal degree  $d$, for 
which  ${u(G)\over v(G)}$ tends to $d-2\over d+1$.

For  $\delta(G)=3$ and $\delta(G)=4$ the statement of Linial's conjecture was proved by   Kleitman and West (\cite{KW}, 1991), for $\delta(G)=5$ --- by Griggs and Wu (\cite{JGM}, 1996). 
In both papers the proof is based on the method of  {\it dead vertices}. 
There are  significant problems with the extension of this method for ${\delta(G)\ge 6}$ and no further results. It follows from the works~\cite{Alon,DJ,YW} that for $\delta(G)$ large enough the Linial's conjecture fails.  However, for small  $\delta(G)>5$ the question remains open.

A number of works consider spanning trees in classes of graphs with various additional conditions
like prohibition  of a certain subgraph. Most of such papers study spanning trees in graphs with no subgraph
isomorphic to $K_4^-$ (complete subgraph on four vertices without one edge).

At first (\cite{KWS}, 1989) it was proved that $u(G)\ge{v(G)+4\over 3}$ for a connected cubic graph without
$K_4^-$. Later Bonsma (\cite{B},2008) proved  two interesting bounds for a connected graph
with $\delta(G)\ge3$: \quad $u(G)\ge{v(G)+4\over 3}$ for a graph without triangles (that is, with $g(G)\ge4$)
and $u(G)\ge{2v(G)+12\over 7}$ for a graph without  $K_4^-$.

These results do not answer the question, how to estimate the number of leaves in spanning tree for
a connected graph with vertices of degrees 1 and 2. Recently some papers were published in which
the  vertices of degree 1 and 2 does not influence the estimate. It is proved in~\cite{BZ}
that $u(G)\ge{v_3+4\over3}$ for a  connected graph $G$ with $g(G)\ge4$ and $v_3$ vertices of degree at
least~3 (indeed, in~\cite[theorem 1]{BZ} this bound was formulated and proved for more general class of graphs). 
In~\cite{Gr} for a connected graph $G$ with  $v_3$ vertices of degree 3 and
$v_4$  vertices of degree at least 4 the estimate~$u(G)\ge{2v_4\over5}+{2v_3\over15}$ is proved.

The first theorem of our paper is a natural continuation of these results. In our theorem  pendant vertices of  considered graph $G$ make their contribution into $u(G)$.

\begin{thm}
\label{t31}
Let $G$ be a connected  graph with  $v(G)\ge 2$ and $s$ vertices of degree not $2$. 
Then there exists a spanning tree of the graph ~$G$  with at least  ${1\over 4}(s-2)+2$  leaves.
\end{thm}

As an easy consequence of this theorem we get a bound of the number of leaves in  spanning trees of graphs without vertices of degree 2.

\begin{cor}
\label{s32}
Let $G$ be a connected  graph  with  $v(G)\ge 2$ and  without vertices of degree $2$. 
Then there is a spanning tree of the graph ~$G$  with at least   ${1\over 4}(v(G)-2)+2$   leaves.
\end{cor}

\begin{defin}
Let us denote by $\ell(G)$ the number of vertices in maximal chain of successively adjacent vertices
of degree 2 in a graph  $G$.
\end{defin}

It was proved in~\cite{K} for a connected  graph $G$ with $\ell(G)\le k$ (where $k\ge 1$), that $u(G)> {1\over 2k+4}v(G)$.  However, the proof from~\cite{K} significantly used, that  $k\ge 1$. The corollary \ref{s32}
shows that for $k = 0$ one can obtain the same bound. We extend the result of~\cite{K} and obtain  new series of bounds, connecting the number of leaves in a spanning tree with the girth of a graph.

\begin{thm}
\label{t3g}
Let  $G$ be a connected graph, $v(G)\ge 2$, $g(G)\ge g$ and $\ell(G)\le k$ (where  $k\ge 1$ is an integer). 
Then there is a spanning tree of the graph ~$G$  with at least   ${\alpha_{g,k}(v(G)-k-2)+2}$ 
 leaves, where
$$\alpha_{g,k}=\left\{ \begin{array}{ll}   {n\over n(k+3)+1} \mbox{ $($where } n= [{g+1\over2}]),  & 
                                              \mbox{ as }  k<g-2, \\
                                           {g-2\over (g-1)(k+2)}, & \mbox{ as } k\ge g-2 \\
                        \end{array} \right. . $$

\end{thm}

The question about such bounds for the case $k=0$ (i.e., for graphs without vertices of degree 2) remains open.

\section{The basic instrument}
In this section we shall formulate and prove a lemma, which will be our main instrument in all proofs. 
This lemma  helps us to reduce the proof of bounds in theorems 1 and 2  when a graph we consider
have cutpoints. This lemma also allows us to construct extremal examples for our bounds from small pieces.

Our method uses the theory of blocks and cutpoints. We recall here the basic definitions of this theory.
For details see~\cite{X} and other books.

\begin{defin}
A {\it cutpoint} of a connected graph~$G$ is  its vertex,  deleting of which makes graph disconnected.

A nonempty graph is {\it biconnected} if it has no cutpoints.
A {\it block} of a graph~$G$ is its maximal (with respect to inclusion) biconnected subgraph.  

A {\it bridge} of a graph~$G$ is its edge, deleting of which makes graph disconnected 
(i.e. an edge, not contained in any cycle).
\end{defin}

\begin{figure}[!hb]
	\centering
		\includegraphics[width=\columnwidth, keepaspectratio]{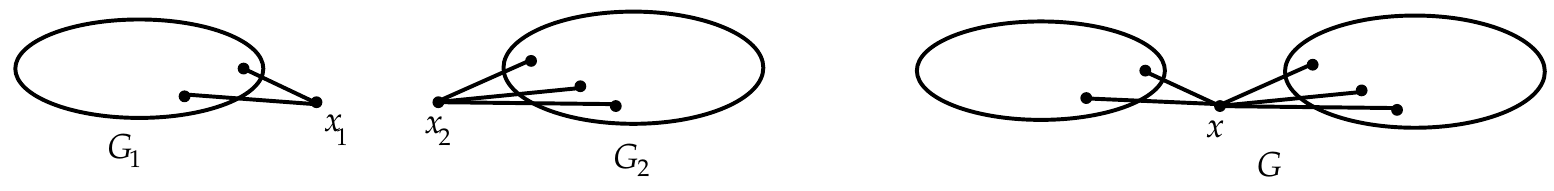}
     \caption{Gluing of graphs.}
	\label{fig1a}
\end{figure}

\begin{defin} $1)$  Let  $G_1$ and $G_2$ be two graphs with marked vertices $x_1\in V(G_1)$ and
$x_2\in V(G_2)$ respectively, $V(G_1)\cap V(G_2)=\varnothing$. 
To {\it glue} the graphs $G_1$ and $G_2$ by the vertices  $x_1$ and $x_2$ is to glue together vertices $x_1$ and $x_2$ into one vertex  $x$,  which will be incident to all edges, incident to  $x_1$ or $x_2$ in both graphs 
$G_1$ and $G_2$. All other vertices and edges of the graphs  $G_1$ and $G_2$  will be vertices and edges of the resulting graph (see fig.\,\ref{fig1a}). 

$2)$  For any edge $e\in E(G)$ we denote by $G\cdot e$ the graph, in which the ends of the edge $e=xy$ are contracted into one vertex, which will be incident to all edges, incident to  $x$ or $y$ in the graph~$G$. Let us say that the graph $G\cdot e$ is obtained from $G$ by {\it contracting} the edge $e$.
\end{defin}

\begin{rem}
\label{rst}
1) Loops and multiple edges do not appear after contracting a bridge.

2) Let a graph  $H$ is obtained from a graph  $H'$ by contracting several bridges, not incident to pendant vertices. Then, obviously, $u(H)=u(H')$.
\end{rem}

\begin{lem}
\label{tool} Let $G_1$ and $G_2$ be connected graphs with  $V(G_1)\cap V(G_2)=\varnothing$, $v(G_1)\ge 2$, 
$v(G_2)\ge 2$  and pendant vertices $x_1$ and $x_2$, respectively. Let  $G$ be a graph, obtained  by gluing 
$G_1$ and $G_2$ by the vertices  $x_1$ and $x_2$  and, after that, by contracting  $m'-1$ bridges, not incident to pendant vertices.  Then the following statements hold.

$1)$ $u(G)=u(G_1)+u(G_2)-2$.

$2)$ Let  $$u(G_1)\ge \alpha(v(G_1)-m)+2, \quad u(G_2)\ge \alpha(v(G_2)-m)+2 \quad 
\mbox{and } m'\ge m.\eqno(1)$$ 
Then $u(G)\ge \alpha(v(G)-m)+2$. If all three inequalities in $(1)$ become equalities, then 
 $u(G)= \alpha(v(G)-m)+2$.

\begin{proof}
1) Let  $G'$ be a graph, obtained  by gluing $G_1$ and $G_2$ by the vertices  $x_1$ and $x_2$.
By remark~\ref{rst} we have $u(G')=u(G)$. It remains to proof, that $u(G')=u(G_1)+u(G_2)-2$.  
Let $x$ be a vertex of the graph $G'$, obtained from $x_1$ and $x_2$ by  gluing together.

$\ge$. Consider spanning trees  $T_1$ and $T_2$ of the graphs  $G_1$ and $G_2$ with  $u(T_1)=u(G_1)$ and
$u(T_2)=u(G_2)$. Gluing together leaves  $x_1$ and $x_2$ into one vertex  $x$, we obtain a spanning tree~$T$ of the graph~$G'$ with $u(T)=u(T_1)+u(T_2)-2$ (all leaves of the trees $T_1$ and $T_2$, except $x_1$ and $x_2$, remain in the tree~$T$). Hence, $u(G')\ge u(G_1)+u(G_2)-2$.

$\le$. Now consider a spanning tree $T'$ of the graph $G'$ with $u(T')=u(G')$.  The vertex~$x$ is a cutpoint of the graph $G'$, and hence $x$ is not a leaf of $T'$, that is, $d_{T'}(x)=d_{G'}(x)=2$. 
It is easy to see, that there exist such spanning trees $T'_1$ of the graph $G_1$ ($x_1$ is a leaf of $T'_1$)
and $T'_2$ of the graph $G_2$ ($x_2$ is a leaf of $T'_2$),  that
the tree  $T'$ is a result of gluing   spanning trees $T'_1$ and $T'_2$ by their leaves   $x_1$ and $x_2$.  All other leaves of the trees  $T'_1$ and $T'_2$  are leaves of the tree $T'$, and thus 
$$u(G')=u(T')=u(T_1')+u(T_2')-2\le  u(G_1)+u(G_2)-2.$$

2) Note, that $v(G)= v(G_1)+v(G_2)-m'$. Indeed, we glue together two vertices $x_1$ and $x_2$ into one 
vertex $x$, after that we contract  $m'-1$ bridges and decrease the number of vertices by  $m'-1$.
After the proof of item  1, it remains only to write  a chain of ineqaulities:
$$u(G)=u(G_1)+u(G_2)-2\ge \alpha(v(G_1)-m)+2 + \alpha(v(G_2)-m)+2 -2 =$$
$$ = \alpha( (v(G_1)+v(G_2)-m')-m +(m'-m)) +2  \ge   \alpha(v(G)-m) +2. \eqno(2)$$
It is easy to see, that in the case when all the inequalities in  $(1)$ become equalities, all the 
inequalities in   $(2)$ become equalities too.
\end{proof}
\end{lem}

\section{Theorem  1 and  series of extremal examples}

\renewcommand*{\proofname}{\bf Proof of Theorem 1}
\begin{proof}
\renewcommand*{\proofname}{\bf Proof}
We  assume that our graph  $G$ has pendant vertices, otherwise we can use the following result of Kleitman and West.

\begin{lem}
\label{kw} {\rm \cite{KW}}
Let $G$ be a connected graph, $\delta(G)\ge 3$. Then $u(G)\ge {v(G) \over 4} +2$.
\end{lem}

Set the following notations:

 $U$ is the set of all pendant vertices of the graph $G$;

 $W$ is the set of all vertices of the graph $G$, which are adjacent to pendant vertices;

 $X$ is the set of all vertices of the graph $G$ not from $U\cup W$, which are adjacent to vertices of $W$;

 $Y$ is the set of all other vertices of the graph $G$. 

For any graph $F$ we denote by $S(F)$ the set of all vertices of this graph of degree not 2. Denote by  $s(F)$
the number of vertices in $S(F)$. Let  $H=G-U$. Obviously, the graph $H$ is  connected.

We  proof  theorem~\ref{t31} by  induction on the size of a graph. 
We assume, that the statement of our theorem is proved for all graphs with less vertices, than $G$, or with the same number of vertices and less number of edges.

\smallskip
{\bf Base  of induction. } The case $v(H)\le 2$ is obvious. It is not difficult to look over all the variants in this case.

\smallskip

{\bf Induction step.}  We consider several cases. 

\smallskip
\q1. {\it There is a vertex   $a\in V(G)$  of degree $2$.}

\noindent 
 If  $a$ is a cutpoint, then, by contracting one edge incident to $a$, we obtain a smaller graph $G'$ with $s(G')=s(G)$ и and $u(G')=u(G)$.

Suppose, that  $a$ is not a cutpoint. Consider an incident to $a$ edge $ab$. This edge is not a bridge, 
thus, the graph $G'= G-ab$ is connected.  Obviously, $a\in S(G')\setminus S(G)$.
Note, that  $S(G)\setminus S(G') \subseteq \{b\}$, since degrees of every different from  $a$ and $b$ vertex in
the graphs $G$ and $G'$ coincide. Hence, $s(G')\ge s(G)$. Since any spanning tree of the graph $G'$ is a spanning tree of the graph $G$, we have $u(G')\le u(G)$. 

In both cases the statement of our theorem for the graph $G$ follows from the statement for a smaller
graph $G'$.

{\it Further  we suppose, that  $s(G)=v(G)$ (i.e. all vertices of the graph  $G$ have degree not $2$).}

\smallskip
\q2. {\it The graph $H$ is not biconnected.}

 Let $a$ be a cutpoint of the graph $H$.    Then $a$ is a cutpoint of the graph $G$, which splits $H$ into 
at least two parts. Thus there exist such connected graphs $G_1$ and $G_2$, that
$$V(G_1)\cup V(G_2)=V(G), \quad V(G_1)\cap V(G_2)=\{a\} \quad \mbox{and} \quad  v(G_1), v(G_2) > 2.$$

For $i\in\{1,2\}$ consider a graph $G_i'$, obtained from $G_i$ by adding a new pendant vertex $x_i$, adjacent to $a$ (see fig.~\ref{fig2a}). We think of two copies of the vertex  $a$ in the graphs $G_1'$ and $G_2'$ as of two different vertices. Then the graph $G$ is obtained from  $G'_1$ and $G_2'$ by gluing together 
the vertices $x_1$ and $x_2$  into one vertex $x$ and contracting two bridges, incident to $x$ (after this process two copies of the vertex  $a$ in graphs $G_1'$ and $G_2'$ are glued into one vertex $a$ of the graph $G$). Thus, the conditions of item  1 of lemma~\ref{tool} are fulfilled and we have $u(G)=u(G_1')+u(G_2')-2$.
It is easy to see, that   $v(G_1')<v(G)$ and $v(G_2')<v(G)$. Then by induction we have $u(G'_1)\ge {s(G'_1)-2\over 4}+2$ and  $u(G'_2)\ge {s(G'_2)-2\over 4}+2$.

\begin{figure}[!hb]
	\centering
		\includegraphics[width=\columnwidth, keepaspectratio]{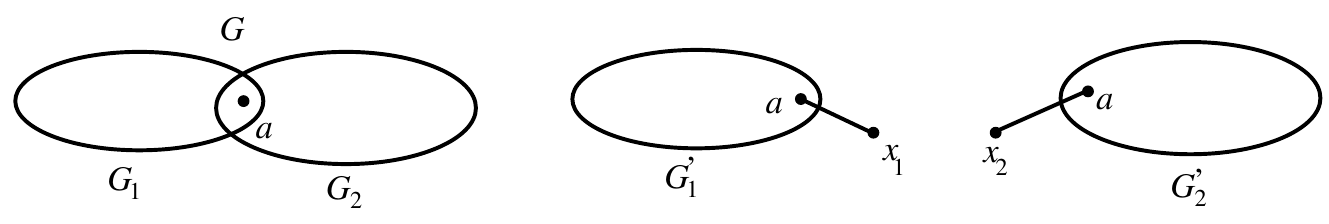}
     \caption{Splitting of the graph by its cutpoint $a$.}
	\label{fig2a}
\end{figure}

Note, that $s(G)=v(G)$ and all vertices of the graphs $G'_1$ and $G'_2$, except $a$, have degree not 2. Since $3\le  d_{G}(a)= d_{G'_1}(a) + d_{G'_2}(a)-2$, the vertex $a$ has degree not 2 in at least one of the graphs
$G_1'$ and $G_2'$. Hence 
$$s(G)=v(G)=v(G_1')+v(G_2')-3\le s(G_1')+s(G_2')-2.$$

Thus, 
$$u(G)=u(G_1')+u(G_2')-2 \ge {s(G'_1)-2\over 4} +  {s(G'_2)-2\over 4} +2  \ge {s(G)-2\over 4}+2,$$
which was to be proved.

\smallskip
{\it In what follows we suppose, that the graph $H$ is biconnected. Hence, all cutpoints of the graph $G$ 
are  vertices of the set  $W$ (any vertex $w\in W$ separates pendant vertices, adjacent to $w$, from all other vertices of the graph~$G$).}

For further progress we need the following lemma.

\begin{lem}
\label{tr1}
Let $a,b\in V(G)$ be adjacent vertices, and subgraph $G'$ be a connected component
of the graph $G-a$, which contains the vertex $b$. Let  $b$ be a cutpoint of the graph $G'$.
Then  $u(G)\ge u(G')+1$.

\begin{proof} Consider a spanning tree  $T'$ of the graph  $G'$ with $u(T')=u(G')$.
Let us construct a spanning tree  $T$ of the graph $G$. We adjoin the vertex  $a$ to  $b$ 
and after that adjoin to the vertex  $a$ all other connected components of the graph $G-a$. The vertex  $b$ is a cutpoint of the graph $G'$, thus, $b$ is not a leaf of the spanning tree $T'$. Hence, $u(G)\ge u(T)\ge u(T')+1=u(G')+1$.
\end{proof}
\end{lem}

Let us continue our analysys of cases in the proof of  theorem~\ref{t31}.

\smallskip
\q3. {\it There exist adjacent vertices $a,b\in V(G)$, such that ${d_G(a)\le 3}$ and $b$ is a cutpoint of the graph  $G-a$.}

\noindent
Let  $G'$ be a connected component of the graph $G-a$, which contains the vertex $b$.
Since the graph  $H$ is biconnected, then  $G'$ contains all vertices of the graph $G$, except
$a$ and pendant vertices, adjacent to $a$. Thus,
$$S(G)\setminus S(G') \subseteq \{a\}\cup \N_G(a), \quad \mbox{hence} \quad
s(G')\ge s(G)-d_G(a)-1\ge s(G)-4.$$
By induction, by  lemma~\ref{tr1} and proved above we have
$$u(G)\ge u(G')+1 \ge  {s(G') -2\over 4}+3 \ge  {s(G) -2\over 4}+2,$$ 
which was to be proved.

\smallskip
\q4. {\it There exist adjacent vertices $x,y\in V(H)$, such that ${d_G(x)\ge 4}$, $d_G(y)\ge 4$.}

\noindent Consider a graph $G'=G-xy$. Since the graph $H$ is biconnected, then the graph $G'$ is connected. The statement of our theorem for smaller graph $G'$ is proved. It is clear, that  $s(G')=s(G)$ and any spanning tree of the graph  $G'$ spanning tree of the graph  $G$, thus the statement is proved for $G$.

\smallskip
\q5.  Summarize the analyzed cases and clear up the properties,  which have the graph $G$, not satisfying the conditions of previous cases.

\begin{lem}
\label{prop}
Let the graph   $G$   does not satisfy the condition of any analyzed  case. Then $G$ has the following properties.

 $1^\circ$   No two  vertices of the set  $W$ are adjacent.

 $2^\circ$    All vertices of the set $W$ have degree $3$.

 $3^\circ$    The set  $X$ is nonempty and consists of vertices of degree more, than~$3$.

 $4^\circ$    Every vertex of the set  $W$ is adjacent to  one pendant vertex of the graph  $G$ and two vertices
of the set $X$.

\begin{proof}
Consider a vertex $w\in W$ and adjacent to $w$ pendant vertex  $u\in U$.
The vertex  $w$ is a cutpoint, separating  $u$ from  other vertices of the graph $G$. 
If the vertex $w$ is adjacent to a vertex of degree not more than 3, different from~$u$, the graph $G$ satisfy the condition of case 3. Thus all vertices, adjacent to  $w$, except the only pendant vertex $u$, have degree more than 3. Hence,  all vertices of the set  $X$ have degree more than 3.

Let us prove, that $W$ is an independent set of the graph  $G$. Let vertices  $w,w'\in W$ be adjacent. Then at least one of them has degree not more than 3, let $d_G(w')\le 3$. But in this case the vertex $w\in W$ is adjacent to non-pendant vertex of degree not more than 3, a contradiction to proved above.

Hence,  $W$ is an independent set.  Then  every vertex   $w\in  W$ is adjacent to  $d_G(w)-1\ge 2$  vertices of the set $X$, thus, the set  $X$ is nonempty.  Since degrees of vertices of $X$ are more than 3, then $d_G(w)=3$. 
All the statements of lemma are proved.
\end{proof}
\end{lem}

 Consider any vertex $w\in W$  and adjacent to $w$ vertices  $x,x'\in X$.
Let  $a\ne w$ be a vertex, adjacent to $x$. It is easy to see, that either $a\in W$, or $a\in Y$ and in both
cases  $d_G(a)=3$ (see fig.~\ref{fig3a}).

\begin{figure}[!hb]
	\centering
		\includegraphics[width=0.4\columnwidth, keepaspectratio]{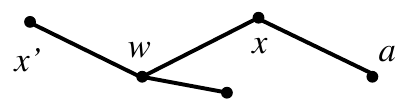}
     \caption{Вершины $w$, $x$, $x'$ и $a$. }
	\label{fig3a}
\end{figure}

Let $G^*=G-x'w$ and  $G'$ be a connected component of the graph  $G^*-a$, which contains $w$. Obviously, the vertex~$x$ is a cutpoint of the graph $G'$ (separating $w$ and a pendant vertex, adjacent to $w$, from another vertices of the graph). Therefore, we can apply  lemma~\ref{tr1} to graphs~$G^*$ and $G'$ and obtain
$u(G^*)\ge u(G')+1$.

The vertex $x$ is adjacent to $a$ and it is not a cutpoint of the graph $G-a$ (otherwise, our graph would be considered in case 3).  Hence, the edge  $xw$ is not a bridge of the graph $G-a$. Since  $w$ is adjacent to
vertices $x,x'\in X$ and one pendant vertex, then $x'w$ also is not a bridge of the graph~${G-a}$. 
Consequently, almost all vertices of the graph $G-a-x'w=G^*-a$ are in the same connected component $G'$
(except  $a$ and --- in the case  $a\in W$ --- a pendant vertex, adjacent to $a$).  

All vertices from  $V(G)$, except  $a,x',w$ and vertices of $\N_G(a)$  are vertices of the graph  $G'$ and 
their degrees in both graphs $G$ and $G'$ coincide. For the vertex  $x\in \N_G(a)$ we have 
$d_{G'}(x)=d_{G}(x)-1\ge3$. If $x'\not\in \N_G(a)$, then $d_{G'}(x')=d_{G}(x')-1\ge3$. Hence,
the set  $S(G)\setminus S(G')$ consists of not more than  $d_G(a)+1\le 4$ vertices: it can be  $w$, $a$ and vertices of $\N_G(a)$, different from $x$.
Hence,  $s(G')\ge  s(G)-4.$
 By induction, $u(G')\ge {s(G') -2\over 4}+2$.
Taking into account, that $G^*$ is a subgraph of $G$ and inequalities, proved above, we have
$$u(G)\ge u(G^*)\ge u(G')+1 \ge  {s(G') -2\over 4}+3 \ge  {s(G) -2\over 4}+2,$$ which was to be proved.
\end{proof}

\subsection{Extremal examples}

\begin{figure}[!hb]
	\centering
		\includegraphics[width=0.8\columnwidth, keepaspectratio]{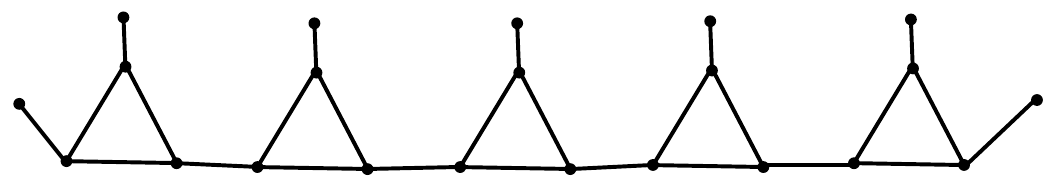}
     \caption{Extremal example for the bound of Theorem 1. }
	\label{fig4a}
\end{figure}

Consider a   tree $T$,  all  vertices of which have  degrees  1 and 3, and exactly $n$ vertices have degree 3.
It is easy to see, that  $T$ has exactly $n+2$ leaves and  $e(T)=2n+1$. 
Substitute every vertex $x$ of degree 3  in the tree $T$ by a triangle: 
each of three vertices of this triangle will be incident to one of  three edges, incident to $x$ in the tree $T$.
An example of such graph for $n=5$ is on fig.~\ref{fig4a}.
The obtained graph $G$ has $n+2$ vertices of degree 1 and  $n$ triangles, altogether $v(G)=n+2+3n=4n+2$. 
All non-pendant vertices of the graph $G$ are its cutpoints, thus they cannot be leaves in a spanning tree. Hence,   $u(G) = n+2= {v(G)-2\over 4}+2$.

\section{Theorem 2 and series of extremal examples}

We begin with necessary definitions and statement of lemma about splitting of large blocks from the paper~\cite{K}.

\begin{defin} Let $B$ be a block of a graph $G$.

The {\it boundary} of    $B$ is the set of all contained in $B$ cutpoints of the graph $G$ (notation: $\R(B)$). The {\it interior} of $B$ is the set of vertices $\I(B)=V(B)\setminus \R(B)$. We call vertices of $\I(B)$ {\it inner vertices} of the block~$B$.

A block $B$ is called   {\it empty}, if it has no inner vertices (i.e., $\I(B)=\varnothing$) . Otherwise, a block is called {\it nonempty}.

A block $B$ is called {\it large}, if the number of its inner vertices is more, than the number of its boundary vertices (i.e., $|\I(B)|>|\R(B)|$).
\end{defin}

\begin{lem}
\label{lrb2}
Let  $G$ be a graph with more than $2$ vertices. Then there exists a set of edges $F\subset E(G)$, satisfying the following conditions:

$1^\circ$  the graph $G-F$ is connected;

$2^\circ$  the graph $G-F$ has no large blocks;

$3^\circ$  if vertices $x$ and $y$ are adjacent in   $G-F$  and $d_{G-F}(x)=d_{G-F}(y)= 2$, then
 $d_{G}(x)=d_{G}(y)= 2$.
\end{lem}

\renewcommand*{\proofname}{\bf Proof of Theorem 2}
\begin{proof}
\renewcommand*{\proofname}{\bf Proof}
The   proof of this  theorem is a descent  on the size of a graph.

\smallskip
\q1. {\it Descent.}

\noindent
We shall say, that a graph $G'$ is {\it smaller} than $G$, if  either $u(G') < u(G)$, or $u(G') = u(G)$ 
and~$e(G') < e(G)$. 
In the first part of proof we analyze  cases, when the statement of Theorem~\ref{t3g} for our graph $G$ follows from the statement for all smaller graphs.

Let a  {\it spine} is a tree with no vertices of degree more than 2, adjoined by an edge from one of its leaves to the cutpoint $a$. The cutpoint $a$ is the {\it base} of this spine.

We call a cutpoint $a$ of the graph  $G$ {\it inessential}, if the graph  $G-a$ has exactly two connected components, and one of these components is a spine with base~$a$. Otherwise we call $a$ an   {\it essential}
cutpoint.

\smallskip
\q{1.1}. {\it A graph $G$ has an essential cutpoint $a$.}

\noindent  
If  $d_G(a)=2$, then the vertex $a$ belongs to some chain of successively adjacent vertices of degree 2, let
 boundary vertices of this chain are adjacent to vertices $b$ and  $b'$ of degree not 2. Since  $a$ is an 
essential cutpoint, then $d_G(b)>2$ and $d_G(b')>2$, vertices $b$ and $b'$ are also essential cutpoints.

Hence it is enough to consider the case  $d_G(a)\ge 3$.  The vertex $a$ is an essential cutpoint of the graph $G$, thus there exist such connected graphs $G_1$ and $G_2$, that $V(G_1)\cup V(G_2)=V(G)$ and $V(G_1)\cap V(G_2)=\{a\}$, moreover, both graphs $G_1$ and $G_2$ are not spines with base $a$. 

Let us construct a graph $G'_1$ from the graph $G_1$. If $d_{G_1}(a) = 1$, then $G_1' = G_1$. If
$d_{G_1}(a) \ge 2$, then adjoin to the vertex  $a$ a spine of $k+1$ vertices. Thus $\ell(G'_1)\le k$, $g(G'_1)\ge g(G)$.  We construct a graph $G'_2$ similarly. 

Since $3\le d_G(a)= d_{G_1}(a)+d_{G_2}(a),$  then $d_{G_1}(a)\ge 2$ or $d_{G_2}(a)\ge 2$. Hence, 
during the construction of at least one of graphs   $G'_1$ or $G'_2$ we add a spine with  $k+1$ vertices. Since 
$a$ is a vertex of both graphs $G'_1$ and $G'_2$, we obtain an ineqaulity $$v(G'_1)+v(G'_2)\ge v(G)+k+2.$$ 
 
The graph $G$ is the result of gluing graphs $G_1'$ and $G_2'$ 
 by two pendant vertices and contracting at least  $k+1$ bridges (after that two copies of the vertex  $a$ in the graphs $G_1'$ and $G_2'$ are glued together into the vertex $a$ of the graph~$G$). By item 1 of lemma~\ref{tool} we have  $u(G)=u(G_1')+u(G_2')-2$. Since the graphs $G_1$ and $G_2$ are not spines with base $a$, then $u(G'_1), u(G'_2) \ge 3$ and, consequently, $u(G_1')<u(G)$ and  $u(G_2')<u(G)$. 
Then we have 
 $$u(G'_1)\ge \alpha_{g,k}(v(G'_1)-k-2)+2,\quad u(G'_2)\ge \alpha_{g,k}(v(G'_2)-k-2)+2.$$
By item  2 of lemma~\ref{tool} we obtain the inequality $u(G)\ge   \alpha_{g,k}(v(G)-k-2)+2,$
which was to be proved.

\smallskip
\q{1.2}. {\it A graph   $G$ has large blocks.}

\noindent
By lemma~\ref{lrb2} we can choose a set of edges $F\subset E(G)$, such that the graph $G'=G-F$ is connected, has no large blocks and for any two vertices $x$ and $y$, adjacent in    $G'$ with  $d_{G'}(x)=d_{G'}(y)= 2$ 
we have  $d_{G}(x)=d_{G}(y)= 2$. Hence $\ell(G')=\max(\ell(G),1)\le k$.  
Obviously, $g(G')\ge g(G)=g$. Thus we can apply the statement of the theorem to smaller graph $G'$. Since any spanning tree of the graph  $G'$ is a spanning tree of the graph $G$, we obtain
$u(G)\ge u(G')\ge \alpha_{g,k}(v(G)-k-2)+2,$ which was to be proved.

\smallskip
\q2. {\it Base.}

\noindent Let us reduce our graph, performing  steps   1.1 and 1.2, until it is possible. It remains to verify the statement of our theorem only for graphs $G$ without essential cutpoints and large blocks. Every cutpoint $a$ of such a graph $G$ splits  $G$ into two connected components and one of these components is a spine with base $a$. Let $H$ be a graph, obtained from $G$ as a result of deleting of vertices of all spines. It is easy to see, that the graph $H$ is biconnected (any cutpoint of $H$ would be an essential cutpoint of the graph $G$).

Let  $h=v(H)$ and $m$ be a number of cutpoints of the graph  $G$. 
Since   $H$ is not a large block of  $G$, then  $m\ge {v(H)\over 2}$. Every cutpoint separates from the graph  a spine with not more than  $\ell(G)+1\le k+1$ vertices. Hence,  $v(G)\le h+(k+1)m$.

The case when  $G$ is a tree we consider individually. It is easy to check, thah all  bounds of theorem~\ref{t3g}
are correct in this case.

Now let $G$ be not a tree, then  biconnected graph  $H$ containes a cycle of not more than $v(H)$ vertices. Consequently, $h=v(H)\ge g(G)=g$. Consider two cases.

\smallskip
\q{2.1}. { $m=h$.}

\noindent Then $v(G)\le (k+2)h$, $u(G)\ge h$. By straightforward calculation it can be checked, that in this case
$$u(G)\ge \beta_{h,k}(v(G)-k-2)+2 \quad \mbox{for} \quad \beta_{h,k}={h-2\over (h-1)(k+2)}.$$
Obviously, $\beta_{h,k}$ increases as  $h$ increases, thus minimal value of $\beta_{h,k}$ is attained at $h=g$ and is equal to $\beta_{g,k}$.
It remains to check, that $\beta_{g,k}\ge \alpha_{g,k},$
which will be done after the case 2.2.

\smallskip
\q{2.2}. { $m<h$.}

\noindent  In this case the block  $H$ is nonempty. Choose a vertex $u\in \I(H)$. It is easy to select in the graph $G$ a spanning tree, which leaves would be  $m$ ends of spines and the vertex $u$, thus  $u(G)\ge m+1$. By straightforward calculation it can be checked, that
$$u(G)\ge \gamma_{h,m,k}(v(G)-k-2)+2 \quad \mbox{for} \quad \gamma_{h,m,k}={m-1\over h+ (k+1)m-k-2}.$$
It remains to check, that for $h\ge g$
$$\beta'_{h,k}=\mmin_{m} \gamma_{h,m,k}\ge \alpha_{g,k}.$$
Note, that $\gamma_{h,m,k}$ increases as $m$ increases, thus minimal value
is attained at  $m=\lceil{h\over2}\rceil$. We obtain, that
 $$\beta'_{h,k}={\lceil{h\over2}\rceil-1\over h+ (k+1)\lceil{h\over2}\rceil-k-2}. $$
It is easy to see, that  $\beta'_{2n-1,k}>\beta'_{2n,k}$.
Straightforward calculation shows us, that
$$\beta'_{2n,k}={n-1\over (n-1)(k+3)+1}$$
increases as  $n$ increases. Thus minimal value of $\beta'_{h,k}$
is attained at $h=2\lceil{g\over2}\rceil$.

\smallskip
{\it Now we must compare the obtained bounds with $\alpha_{g,k}$. Instead of this we compare the bounds,
obtained in items $2.1$ and $2.2$, with each other and show, that $\alpha_{g,k}$ is equal to minimal of them.} We consider the cases of even and odd~$g$ separately.

\smallskip
\q{2a}. { $g=2n+2$.}

\noindent In this case we need to compare $\beta'_{2n+2,k}$ and $\beta_{2n+2,k}$.
Simple calculation shows, that 
$$\beta'_{2n+2,k}={n\over n(k+3)+1}<{2n\over (2n+1)(k+2)} =\beta_{2n+2,k}$$ if and only if  $k< 2n$.
Exactly in the case  $k< 2n=g-2$ we have $\alpha_{g,k}=\beta'_{2n+2,k}$, and in the case $k\ge g-2$ we have $\alpha_{g,k}=\beta_{g,k}$. Therefore, for even $g$ the teorem is proved.

\smallskip \goodbreak
\q{2b}. { $g=2n+1$.}

\noindent In this case we need to compare $\beta'_{2n+2,k}$ and $\beta_{2n+1,k}$.
Simple calculation shows, that 
$$\beta'_{2n+2,k}={n\over n(k+3)+1}< {2n-1\over  2n(k+2)} =\beta_{2n+1,k}$$ if and only if $k< 2n-1-{1\over n}$. We remind, that~$k$ is integer.
Exactly in the case  $k< 2n-1=g-2$ we have $\alpha_{g,k}=\beta'_{2n,k}$, and in the case $k\ge g-2$ we have $\alpha_{g,k}=\beta_{g,k}$. Therefore,  analysis of the case of odd  $g$ finishs the  proof of the theorem.
\end{proof}

\subsection{Extremal examples}

We describe infinite series of examples, showing that all bounds of theorem~\ref{t3g} are tight. Our reasoning  is rather simple: we construct a graph, for which all inequalities, proved in the theorem become equalities. Let $\ell(G)=k$, $g(G)=g$. Consider two cases.

\q1. {$k< g-2$.}

\begin{figure}[!hb]
	\centering
		\includegraphics[width=\columnwidth, keepaspectratio]{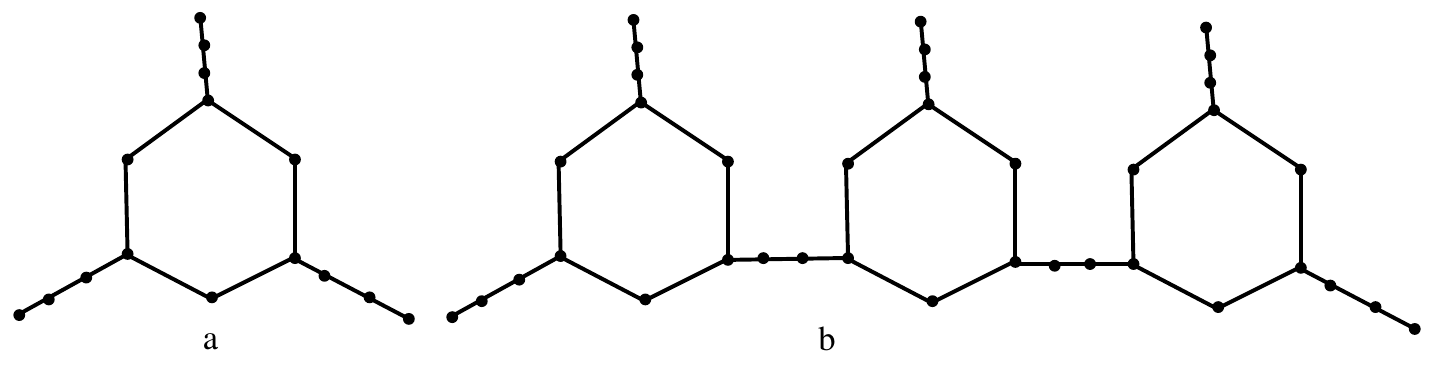}
     \caption{Extremal example for  $k< g-2$. }
	\label{fig5a}
\end{figure}

\noindent Let $n=[{g+1\over2}]$.  In this case $\alpha_{g,k}=\beta'_{2n+2,k}={n\over n(k+3)+1}$.
Let  $B_{g,k}$  be a cycle of length $2n+2$  with  $n+1$ marked vertices (no two marked vertices are adjacent),
and a spine of $k+1$ vertices adjoined to each marked vertex. 
All marked vertices are cutpoints of our graph. Clearly, $\ell(B_{g,k})=k$, $g(B_{g,k})=2n+2 \ge g$.
Then $v(B_{g,k})=2n+2+(n+1)(k+1)=(n+3)(k+1)$. An example of such graph for  $n=2$ (i.e., $g=5$ or $g=6$) and $k=2$ is shown on fig.\,\ref{fig5a}a.

Let us find  $u(B_{g,k})$. All $n+1$ pendant vertices of $B_{g,k}$ (ends of spines) will be leaves of any spanning tree of the graph $B_{g,k}$. Since deleting of leaves of a spanning tree cannot break connectivity of the graph, only one non-pendant vertex of $B_{g,k}$  can be a leaf of its spanning tree (of course, this vertex would be not a base of spine). Thus, $u(B_{g,k})=n+2$. It is easy to check, that
$$u(B_{g,k})=n+2=   2+ {n\over n(k+3)+1} \cdot (v(B_{g,k})-k-2).$$ 
Hence, for the graph  $B_{g,k}$ the bound of theorem~\ref{t3g} is tight.

\smallskip

\q2. {$k\ge  g-2$.}

\begin{figure}[!hb]
	\centering
		\includegraphics[width=\columnwidth, keepaspectratio]{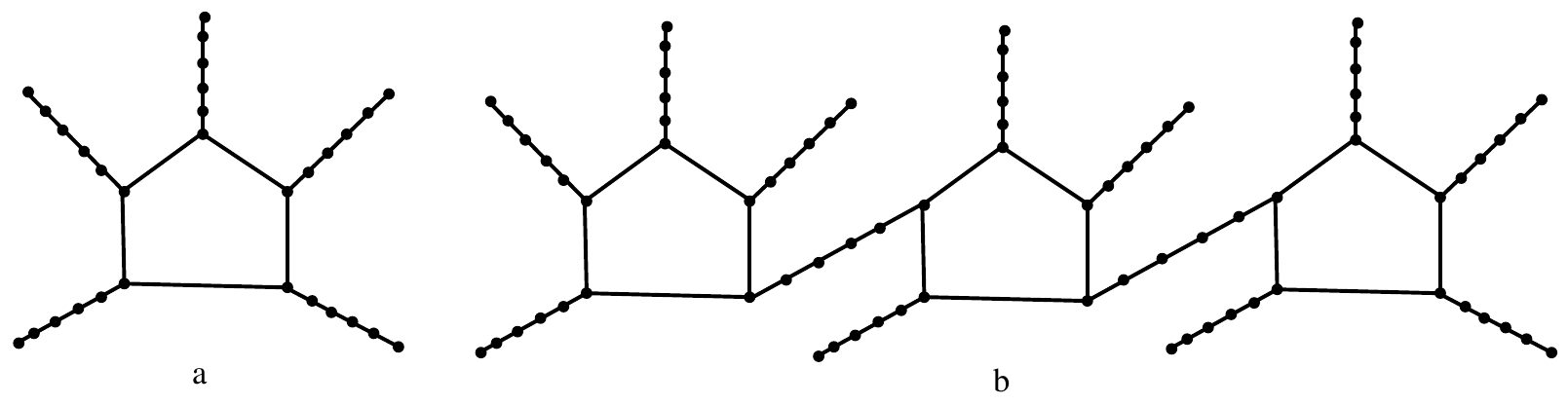}
     \caption{Extremal example for $k\ge g-2$. }
	\label{fig6a}
\end{figure}

\noindent   In this case $\alpha_{g,k}=\beta_{g,k}={g-2\over (g-1)(k+2)}$.
Let  $B_{g,k}$  be a cycle of length~$g$ with  a spine of $k+1$ vertices adjoined to each  vertex. 
Then $v(B_{g,k})=g(k+2)$.   Clearly, $\ell(B_{g,k})=k$, $g(B_{g,k})=g$.
An example of such graph for   $g=5$ and $k=4$ is shown on fig.\,\ref{fig6a}a.

It is easy to see, that ends of spines and only them will be leaves of any spanning tree of this graph, thus $u(B_{g,k})=g$. It is easy to check, that 
$$u(B_{g,k})=g=   2+ {g-2\over (g-1)(k+2)} \cdot (v(B_{g,k})-k-2).$$ 
Hence, for the graph  $B_{g,k}$ the bound of theorem~\ref{t3g} is tight.

\q3. {\it Let us show how to construct extremal examples from the details $B_{g,k}$ in both cases.}

\noindent Let $G$ be a graph, satisfying the conditions
$$u(G)=    \alpha_{g,k} \cdot (v(G)-k-2) +2, \quad g(G)\ge g,\quad \ell(G)\le k,$$ 
and having at least one pendant vertex $a$. We  construct a graph  $G'$:  glue together the vertex $a$ of the graph $G$ with an end of one spines of the graph  $B_{g,k}$ and, after that, contract  $k+1$ bridges  
(edges of glued spine of $B_{g,k}$).  As a result we obtain the graph $G'$, satisfying the following conditions:
$$v(G')=v(G)+v(B_{g,k})-k-2, \quad g(G')\ge g,\quad \ell(G')\le k.$$ 
By item 2 of lemma~\ref{tool} we have $u(G')=    \alpha_{g,k} \cdot (v(G')-k-2) +2,$ i.e., the graph $G'$ is 
also an extremal example, which confirms that the bound of theorem~\ref{t3g} is tight. 
At first we take  $G=B_{g,k}$, then we can construct arbitrary large extremal examples, each time gluing a next graph  $B_{g,k}$ to the graph we have.
Two such examples are shown on fig.~\ref{fig5a}b and~\ref{fig6a}b.

\medskip
Translated by D.\,V.\,Karpov.

\end{document}